\documentclass[review]{elsarticle}
\usepackage[top=3cm,bottom=3cm,left=3cm,right=3cm]{geometry}
\usepackage{amsmath,amssymb,amsthm}
\usepackage{xcolor}
\usepackage{tikz}

\usepackage{array,longtable}
\usepackage{algorithm}
\usepackage{algpseudocode}
\biboptions{numbers,sort&compress}
\usepackage{hyperref}
\definecolor{edgezero}{HTML}{2166AC}
\definecolor{edgeone}{HTML}{D66000}
\definecolor{edgetwo}{HTML}{008060}
\tikzset{av/.style={circle,draw,fill=white,inner sep=2pt,minimum size=7mm},
 bv/.style={rectangle,draw,fill=white,inner sep=2pt,minimum size=7mm},
 ee/.style={line width=1.2pt}, every node/.style={font=\small}}
\newtheorem{theorem}{Theorem}
\newtheorem{lemma}[theorem]{Lemma}
\newtheorem{proposition}[theorem]{Proposition}
\newtheorem{corollary}[theorem]{Corollary}
\newtheorem{conjecture}[theorem]{Conjecture}
\theoremstyle{definition}
\newtheorem{definition}[theorem]{Definition}
\newtheorem{example}[theorem]{Example}
\newtheorem{question}[theorem]{Question}

\theoremstyle{remark}

\begin{document}
\begin{frontmatter}
\title{Properly Colored Hamiltonian Paths \\ in Edge-Colored Complete Bipartite Graphs\tnoteref{funding}}
\tnotetext[funding]{This work is supported by the Scientific and Technological Research Council of Turkey (TUBITAK) under grant no.\ 124F435.}
\author[gtu]{Yasemin B\"uy\"uk\c{c}olak}
\ead{y.buyukcolak@gtu.edu.tr}
\author[gtu]{Fatih Yetgin}
\ead{fyetgin@gtu.edu.tr}
\address[gtu]{Department of Mathematics, Gebze Technical University, Kocaeli, Turkey}

\begin{abstract}
A subgraph of an edge-colored graph is properly colored if no two
adjacent edges have the same color. The Bang-Jensen-Gutin Theorem provides a
characterization of properly colored Hamiltonian paths in complete graphs $K_{n}$
by the existence of a properly colored $1$-path-cycle factor. The same characterization was also conjectured for
complete balanced bipartite graphs $K_{n,n}$. In this paper, we disprove this conjecture
by constructing counterexamples on $K_{n,n}$
for every $n\geq 7$, even with only three colors and with minimum color-degree two.
We then establish the characterization for
colorings in which every {monochromatic color component} is a
path, an even cycle, or a complete bipartite graph with a matching
removed.
The characterization also holds for the complete bipartite graphs $K_{m,n}$ with
$|m-n|=1$. Given a properly colored $1$-path-cycle factor, our
algorithm constructs a properly colored Hamiltonian path in $O(mn)$
time, which is linear in the size of the input.
\end{abstract}
\begin{keyword}
edge-colored graph \sep properly colored Hamiltonian path \sep
complete bipartite graph \sep 1-path-cycle factor \sep
polynomial-time algorithm
\MSC[2020] 05C15 \sep 05C38 \sep 05C45 \sep 05C85
\end{keyword}
\end{frontmatter}

\section{Introduction}\label{sec:introduction}

A subgraph of an edge-colored graph is called \emph{properly colored}
(PC) if no two adjacent edges have the same color. Properly colored
paths and cycles have been extensively studied in edge-colored graph
theory, particularly in connection with spanning structures; see
the surveys~\cite{BJG,GutinKim}. A central problem in this area is
determining when an edge-colored graph contains a properly colored
Hamiltonian path (PCHP) or a properly colored Hamiltonian cycle
(PCHC). Various sufficient conditions for the existence of PCHCs
have been established, notably for complete
graphs~\cite{BollobasErdos,Shearer,AlonGutin,LoHamilton} and, under
color-degree conditions, for general edge-colored
graphs~\cite{LoDirac}.

A \emph{PC $1$-path-cycle factor} is a spanning collection of
pairwise vertex-disjoint PC subgraphs consisting of exactly one path
and zero or more cycles. Clearly, every PCHP is such a factor.
Bang-Jensen and Gutin~\cite{BJG} conjectured that the converse holds
in complete graphs, and this was confirmed in 2006 by Feng et
al.~\cite{Feng}. The resulting characterization, known as the
Bang-Jensen-Gutin (BJG) Theorem, states that an edge-colored
complete graph contains a PCHP if and only if it contains a PC
$1$-path-cycle factor. Together with polynomial-time algorithms for
finding such factors~\cite{BJG}, this characterization yields a
polynomial-time method for deciding the existence of a PCHP in
edge-colored complete graphs.

A natural question is whether this characterization extends to
complete balanced bipartite graphs. In this direction, Cheng et
al.~\cite{CKW}
proposed the following conjecture in 2020.

\begin{conjecture}[Bipartite Extension of the BJG Theorem]
\label{conj:bipartite-BJG}
An edge-colored complete balanced bipartite graph contains
a PCHP if and only if it contains a PC $1$-path-cycle factor.
\end{conjecture}

\noindent
The bipartite setting introduces additional difficulties. Unlike
in complete graphs, paths in bipartite graphs must alternate between
the two partite sets. Consequently, transforming a PC $1$-path-cycle
factor into a single PC spanning path requires suitable connections
between its components that respect both the bipartition and the
proper coloring.

Several partial results toward Conjecture~\ref{conj:bipartite-BJG}
have been obtained in recent years. Büyükçolak~\cite{BYTwo} verified
the conjecture for two-edge-colored complete balanced bipartite
graphs and provided a constructive polynomial-time algorithm.
In another direction, Büyükçolak and Yetgin~\cite{BYC4} proved that every
edge-colored complete balanced bipartite graph containing neither
monochromatic nor semi-monochromatic $4$-cycles admits a PCHP,
without restricting the number of colors.

In this paper, we show that the distinction between two and three
colors is decisive for Conjecture~\ref{conj:bipartite-BJG}. Although
the conjecture holds for two colors~\cite{BYTwo}, it fails already
for three colors: we construct a three-edge-colored $K_{7,7}$ with
minimum color-degree two that admits a PC $1$-path-cycle factor but
contains no PCHP, and we extend the construction to $K_{n,n}$ for
every $n\geq7$ (Theorem~\ref{thm:counterexample}). Consequently,
three is the minimum number of colors for which the conjecture
fails.

Motivated by this counterexample, we seek structural conditions on
the edge-coloring under which the existence of a PC $1$-path-cycle
factor guarantees a PCHP. A natural approach is to examine the color
classes. For a color $i$, consider the subgraph formed by the edges
of color $i$; each of its connected components containing at least
one edge is called a \emph{monochromatic color component} (see
Section~\ref{sec:preliminaries}). We focus on the class $\mathcal F$
of edge-colored bipartite graphs in which every monochromatic color
component is a path, an even cycle, or a complete bipartite graph
with a possibly empty matching removed (for instance, a crown
graph, that is, $K_{t,t}$ minus a perfect matching).

Our main result (Theorem~\ref{thm:component-characterization})
establishes that a complete bipartite graph $K_{m,n}$ in
$\mathcal F$ has a PCHP if and only if it has a PC $1$-path-cycle
factor. This characterization holds not only in the balanced case
$m=n$ of the original conjecture, but also in the nearly balanced
case $|m-n|=1$, where both endpoints of the path lie in the larger
partite set. Since both a PCHP and a PC $1$-path-cycle factor
require $|m-n|\le 1$, these two cases cover all complete bipartite
graphs in which either structure can exist. Moreover, our proof is
constructive: starting from any PC $1$-path-cycle factor, we obtain
a PCHP by successively absorbing its cycles in an arbitrary order.
Algorithmically, we show that membership in $\mathcal F$ can be
recognized in $O(mn)$ time and that a given PC $1$-path-cycle factor
can be transformed into a PCHP within the same time bound. Together
with existing polynomial-time algorithms for finding PC
$1$-path-cycle factors~\cite{BJG}, these results yield a
polynomial-time algorithm for deciding the existence of a PCHP in
$\mathcal F$.

The rest of the paper is organized as follows.
Section~\ref{sec:preliminaries} presents the necessary preliminaries,
and Section~\ref{sec:counterexamples} presents our counterexample to
Conjecture~\ref{conj:bipartite-BJG}.
Section~\ref{sec:monochromatic-components} establishes our main
characterization, while Section~\ref{sec:algorithms} develops the
corresponding algorithms. Concluding remarks are given in
Section~\ref{sec:scope}.

\section{Preliminaries}
\label{sec:preliminaries}

All graphs considered in this paper are finite, simple, and undirected.
For a graph $G$, let $V(G)$ and $E(G)$ denote its vertex and edge sets,
respectively, and let $G[X]$ denote the subgraph induced by
$X\subseteq V(G)$. For a subgraph $H$ of $G$ and $v\in V(H)$, we write
$d_H(v)$ for the degree of $v$ in $H$ and $\Delta(H)$ for the maximum
degree of $H$. The \emph{order} of a path or cycle is its number of
vertices, and its \emph{length} is its number of edges. We write $P_t$
and $C_t$ for a path and a cycle of order $t$, respectively, and denote
a path by the sequence of its vertices, such as $x_0x_1\cdots x_{t-1}$.
A path of order one is called \emph{trivial}. A \emph{segment} of a path
is a subpath formed by consecutive vertices. The first and last edges
of a nontrivial path, with respect to a chosen orientation, are called
its \emph{terminal edges}; for a path of length one, both refer to the
same edge. A \emph{Hamiltonian path} or \emph{Hamiltonian cycle} is a
path or cycle containing every vertex of the graph.
Unless otherwise stated, $G=K_{m,n}$ is a complete bipartite graph
with bipartition $(A,B)$, where $|A|=m\geq1$ and $|B|=n\geq1$.
A vertex set $X\subseteq V(G)$ is \emph{balanced} if
$|X\cap A|=|X\cap B|$.

An \emph{edge-coloring} of $G$ is a map $c:E(G)\to\mathcal C$, where
$\mathcal C$ is a set whose elements are called \emph{colors}; we
usually take $\mathcal C$ to be a set of nonnegative integers. The
coloring need not be proper, that is, adjacent edges may receive the
same color. A graph together with such a map is called an
\emph{edge-colored graph}, and it is \emph{$k$-edge-colored} if exactly
$k$ colors occur on its edges, that is, $|c(E(G))|=k$. For a subgraph
$H$ of $G$, let $c(H)=\{c(e):e\in E(H)\}$ denote the set of colors
appearing on its edges. When $A=\{a_0,\ldots,a_{m-1}\}$ and
$B=\{b_0,\ldots,b_{n-1}\}$, the coloring is represented by the
$m\times n$ \emph{color matrix} $M$ with entries $M_{ij}=c(a_ib_j)$;
thus rows correspond to vertices of $A$ and columns to vertices of $B$.
The \emph{color-degree} $d^c(v)$ of a vertex $v$ is the number of
distinct colors appearing on its incident edges, and the
\emph{minimum color-degree} of $G$ is
$\delta^c(G)=\min_{v\in V(G)}d^c(v)$.

A subgraph of an edge-colored graph is \emph{properly colored} (PC) if
any two edges sharing a vertex have distinct colors, and
\emph{monochromatic} if all its edges have the same color. A $4$-cycle
is \emph{semi-monochromatic} if its edges have cyclic color sequence
$(a,a,b,b)$, where $a\ne b$. A PC Hamiltonian path and a PC Hamiltonian
cycle are abbreviated as PCHP and PCHC, respectively.
For a color $i\in\mathcal C$, let $G_i$ denote the spanning subgraph
of $G$ whose edges are exactly the edges of color $i$; we call $G_i$
the \emph{color class} of $i$. A \emph{monochromatic color component}
of color $i$ is a connected component of $G_i$ containing at least one
edge. Thus monochromatic color components are taken within the color
classes, not within $G$ itself, and a single color class may have
several monochromatic color components.

A \emph{PC $1$-path-cycle subgraph} of $G$ is a collection
$F=\{P,C_1,\ldots,C_k\}$ of pairwise vertex-disjoint PC subgraphs,
where $P$ is a  path, $C_1,\ldots,C_k$ are cycles,
and $k\geq0$. It is a \emph{PC $1$-path-cycle factor} if, in addition,
it spans $V(G)$. Every cycle contains equally many vertices in the two
partite sets. Hence the path component accounts for the difference
$m-n$, and such a factor can exist only if $|m-n|\leq1$. The path has
even order when $m=n$ and odd order when $|m-n|=1$; in the latter case
it may be trivial. When $m=n$, the endpoints of $P$ lie in opposite
partite sets. When $|m-n|=1$, both endpoints, or the unique vertex of a
trivial path, lie in the larger partite set.

For a specified cycle $C$, the notation $E(C)$ refers to its cycle
edges, rather than all edges of $G[V(C)]$. For $e\in E(C)$, the path
$C-e$ is a Hamiltonian path of $G[V(C)]$, called the \emph{cut path} of
$C$ at $e$; if $C$ is PC, so is every cut path. Let $Q$ be a PC path
and $C$ a PC cycle disjoint from it. \emph{Attaching} $C$ to $Q$ means
joining a cut path of $C$ to an endpoint of $Q$ by an edge of $G$, and
\emph{inserting} $C$ into $Q$ means replacing an edge $xy$ of $Q$ by a
cut path of $C$ together with two edges joining its endpoints to $x$
and $y$. The edges added in either operation are called
\emph{joining edges}.

For a vertex set $Y\subseteq V(G)$, a \emph{$Y$-block} of a path is a
maximal nonempty segment whose vertices belong to $Y$. Such a block is
\emph{internal} if it contains neither endpoint of the path, and \emph{terminal} otherwise. The
following elementary observation, related to the notion
of a degenerate set introduced in~\cite{LLZ} (see also~\cite{LiDeg}
and~\cite{LN}), will be useful
in establishing our counterexamples.

\begin{lemma}\label{lem:labels}
Let $Y\subseteq V(G)$ admit a labeling $\ell:Y\to\mathcal C$ such that
$c(yz)\in\{\ell(y),\ell(z)\}$ for every $yz\in E(G[Y])$. Then no PC
path contains a segment $xy_1\cdots y_sz$ with $y_1,\ldots,y_s\in Y$,
$x,z\notin Y$, $c(xy_1)=\ell(y_1)$, and $c(y_sz)=\ell(y_s)$.
In particular, if $X'\subseteq V(G)\setminus Y$ satisfies
$c(xy)=\ell(y)$ for every edge $xy$ with $x\in X'$ and $y\in Y$, then
no PC path contains a $Y$-block with vertices of $X'$ immediately on
both sides.
\end{lemma}

\begin{proof}
Suppose that a PC path contains such a segment.
By assumption, $c(xy_1)=\ell(y_1)$. Since the path is PC and
$c(y_iy_{i+1})\in\{\ell(y_i),\ell(y_{i+1})\}$, it follows inductively
that $c(y_iy_{i+1})=\ell(y_{i+1})$ for $1\leq i<s$. Consequently, the
edge entering $y_s$ has color $\ell(y_s)$, as does the edge $y_sz$,
contradicting proper coloring. For $s=1$, the edges $xy_1$ and $y_1z$
both have color $\ell(y_1)$ directly.
\end{proof}
\section{A Counterexample to the Bipartite Extension Conjecture}
\label{sec:counterexamples}

In this section, we show that the existence of a PC $1$-path-cycle
factor does not guarantee a PCHP in an edge-colored complete balanced
bipartite graph. Although the conjecture holds for two
colors~\cite{BYTwo}, it fails already for three colors, even with
minimum color-degree two. We give an explicit counterexample on
$K_{7,7}$ and extend the construction to $K_{n,n}$ for every $n\geq7$.

\begin{theorem}\label{thm:counterexample}
For every integer $n\geq7$, there exists a three-edge-colored
$K_{n,n}$ with $\delta^c(G)=2$ that admits a PC $1$-path-cycle factor
but has no PCHP.
\end{theorem}

\begin{proof}
We first construct the example for $n=7$. Let $G=K_{7,7}$ have
bipartition $A=\{a_0,\ldots,a_6\}$ and $B=\{b_0,\ldots,b_6\}$, with
edge-coloring given by the color matrix
\[
\begin{array}{c|ccccccc}
 & b_0&b_1&b_2&b_3&b_4&b_5&b_6\\ \hline
a_0&0&0&0&1&0&0&0\\
a_1&1&0&1&1&1&0&0\\
a_2&1&2&1&1&1&0&0\\
a_3&2&0&1&2&1&0&2\\
a_4&0&0&0&1&0&0&0\\
a_5&1&1&1&1&1&0&1\\
a_6&1&1&1&1&1&1&0
\end{array}.
\]
\smallskip

\noindent
The graph admits the PC $1$-path-cycle factor $F=\{P,C_1,C_2\}$, where
\(P=a_0b_0a_1b_1a_2b_2\), \(C_1=a_3b_3a_4b_4a_3\), and
\(C_2=a_5b_5a_6b_6a_5\), with respective color sequences
$(0,1,0,2,1)$, $(2,1,0,1)$, and $(0,1,0,1)$; see
Figure~\ref{fig:counterfactor}. Every row and every column of the
color matrix contains at least two distinct colors, with equality,
for instance, for $a_0$. Thus $\delta^c(G)=2$.

It remains to show that $G$ has no PCHP. Suppose, to the contrary,
that $H$ is a PCHP of $G$. Since $m=n$, the two endpoints of $H$ lie
in opposite partite sets.
Notice here that the rows of $a_0$ and $a_4$ coincide: each of these vertices has a
unique edge of color different from $0$, namely $a_0b_3$ and
$a_4b_3$, both of color $1$. We call $a_0$ and $a_4$ \emph{twins}.
An internal vertex of $H$ is incident with two edges of $H$ of
distinct colors, so an internal twin uses its edge to $b_3$. If both
twins were internal, color $1$ would be repeated at $b_3$. Since the
endpoints of $H$ lie in opposite partite sets, exactly one twin is an
endpoint of $H$. Denote this twin by $t_e$ and the other by $t_i$,
and orient $H$ from $t_e$. In particular, $t_ib_3\in E(H)$, and the
other endpoint of $H$ lies in $B$.

Let $Y=\{a_5,a_6,b_5,b_6\}$ and assign the labels
$\ell(a_5)=\ell(a_6)=1$ and $\ell(b_5)=\ell(b_6)=0$. Every edge of
$G[Y]$ has the label of one of its endpoints. Every edge joining $Y$
to $V(G)\setminus Y$ has the label of its endpoint in $Y$, except
$a_3b_6$, which has color $2$. By Lemma~\ref{lem:labels}, at
least one of the two edges joining an internal $Y$-block of $H$ to the
rest of $H$ is $a_3b_6$. We
distinguish two cases according to the position of $b_3$ on $H$.

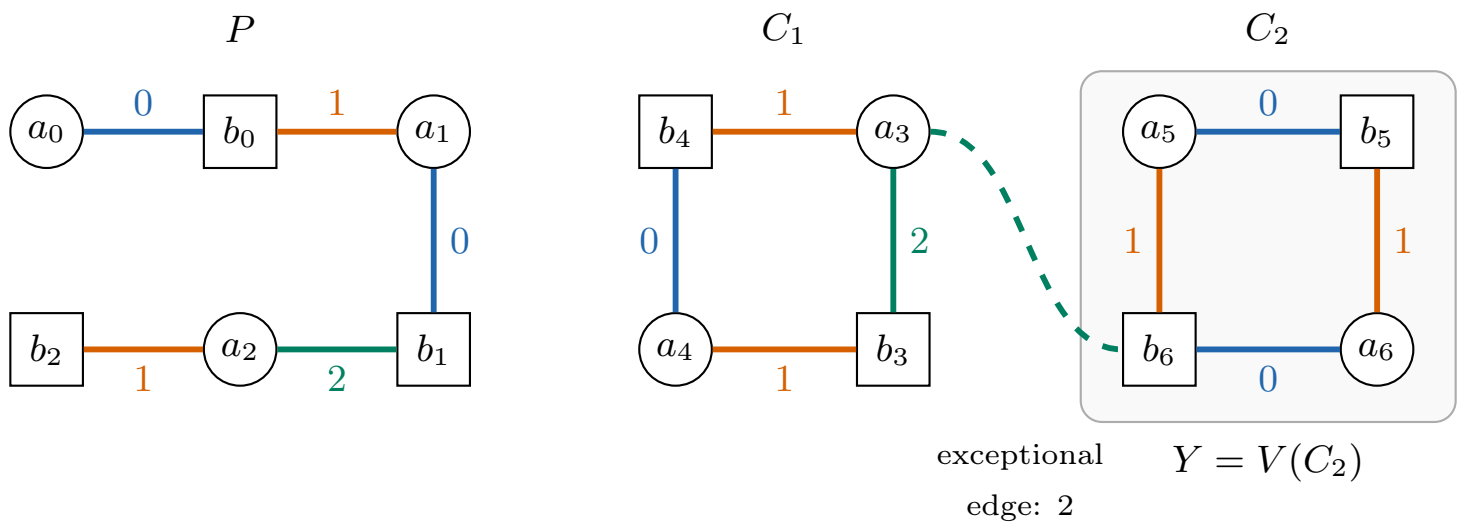
\begin{figure}[tbp]
\centering
\begin{tikzpicture}[
    scale=0.8, transform shape,
    x=1cm,y=1cm,
    every node/.style={font=\small},
    av/.style={circle,draw,fill=white,minimum size=6mm,inner sep=1pt},
    bv/.style={rectangle,draw,fill=white,minimum size=6mm,inner sep=1pt},
    fac/.style={line width=1.1pt}
]
\node[av] (a0) at (0,0) {$a_0$};
\node[bv] (b0) at (1.6,0) {$b_0$};
\node[av] (a1) at (3.2,0) {$a_1$};
\node[bv] (b1) at (3.2,-1.8) {$b_1$};
\node[av] (a2) at (1.6,-1.8) {$a_2$};
\node[bv] (b2) at (0,-1.8) {$b_2$};

\draw[fac,edgezero] (a0)--node[above]{$0$}(b0);
\draw[fac,edgeone] (b0)--node[above]{$1$}(a1);
\draw[fac,edgezero] (a1)--node[right]{$0$}(b1);
\draw[fac,edgetwo] (b1)--node[below]{$2$}(a2);
\draw[fac,edgeone] (a2)--node[below]{$1$}(b2);

\node at (1.6,.85) {$P$};

\node[bv] (b4) at (5.2,0) {$b_4$};
\node[av] (a4) at (5.2,-1.8) {$a_4$};
\node[bv] (b3) at (7.0,-1.8) {$b_3$};
\node[av] (a3) at (7.0,0) {$a_3$};

\draw[fac,edgetwo] (a3)--node[right]{$2$}(b3);
\draw[fac,edgeone] (b3)--node[below]{$1$}(a4);
\draw[fac,edgezero] (a4)--node[left]{$0$}(b4);
\draw[fac,edgeone] (b4)--node[above]{$1$}(a3);

\node at (6.1,.85) {$C_1$};

\draw[draw=gray!65,fill=gray!5,rounded corners]
    (8.55,-2.4) rectangle (11.65,.5);

\node[av] (a5) at (9.2,0) {$a_5$};
\node[bv] (b5) at (11.0,0) {$b_5$};
\node[av] (a6) at (11.0,-1.8) {$a_6$};
\node[bv] (b6) at (9.2,-1.8) {$b_6$};

\draw[fac,edgezero] (a5)--node[above]{$0$}(b5);
\draw[fac,edgeone] (b5)--node[right]{$1$}(a6);
\draw[fac,edgezero] (a6)--node[below]{$0$}(b6);
\draw[fac,edgeone] (b6)--node[left]{$1$}(a5);

\node at (10.1,.85) {$C_2$};
\node[anchor=north] at (10.1,-2.45) {$Y=V(C_2)$};

\draw[fac,edgetwo,dashed]
    (a3.east) .. controls (8.1,0) and (8.1,-1.8) .. (b6.west);
\node[anchor=north,font=\scriptsize,align=center]
    at (8.05,-2.45) {exceptional\\edge: $2$};
\end{tikzpicture}
\caption{A PC $1$-path-cycle factor in the three-edge-colored
$K_{7,7}$ counterexample.}
\label{fig:counterfactor}
\end{figure}

\begin{itemize}
\item \emph{Case 1: $b_3$ is internal.}
Since $c(t_ib_3)=1$, the other edge of $H$ at $b_3$ has a color
different from $1$. In the column of $b_3$, the only such entry is
$c(a_3b_3)=2$, so $b_3a_3\in E(H)$. Hence $a_3b_6\notin E(H)$, since
otherwise color $2$ would be repeated at $a_3$, and therefore $H$ has
no internal $Y$-block. As $t_e\notin Y$, all vertices of $Y$ form a
single terminal $Y$-block containing the endpoint of $H$ in $B$.
Because $Y$ is balanced and this block ends in $B$, it begins at a
vertex of $Y\cap A$.
\item 
\emph{Case 2: $b_3$ is an endpoint.}
Then both endpoints of $H$ lie outside $Y$, so every $Y$-block of $H$
is internal and is joined to the rest of $H$ by $a_3b_6$. Since this
edge can serve only one $Y$-block, all vertices of $Y$ form a single
internal block. Because $Y$ is balanced and
the exceptional edge
$a_3b_6$ joins the block at $b_6$, the other end of the block lies in
$Y\cap A$, and its other joining edge goes from
this vertex to a vertex of $B\setminus Y$.
\end{itemize}

\noindent
In either case, $a_3$ is incident with an edge of $H$ of color $2$,
and no vertex of $A\setminus\{a_3\}$ outside $Y$ is adjacent along $H$
to a vertex of $Y\cap B$. The vertices $a_1$ and $a_2$ are internal,
since the endpoints of $H$ are $t_e$ and a vertex of $B$. Moreover,
they are not adjacent along $H$ to $b_3$, whose neighbors on $H$ are
already $t_i$ and, in Case~1, $a_3$. Their edges to
$W=\{b_0,b_2,b_4\}$ all have color $1$, whereas their only remaining
admissible edges, of colors different from $1$, go to $b_1$.
Therefore $H$ contains the segment $a_1b_1a_2$, up to reversal. Let
$u,v\in W$ be the other neighbors of $a_1$ and $a_2$ on $H$,
respectively. They are distinct, and $c(a_1u)=c(a_2v)=1$. Both $u$
and $v$ are internal, since the endpoint of $H$ in $B$ lies in $Y$ in
Case~1 and equals $b_3$ in Case~2. Hence the other neighbor of each of $u$ and $v$ on $H$ is joined to it
by an edge of color different from $1$. In the columns of $W$, the
only entries different from $1$ are $0$ in the rows of the twins and
$c(a_3b_0)=2$. The latter edge cannot be used, since $a_3$ is already
incident with an edge of color $2$ on $H$. Thus $u$ and $v$ are
adjacent along $H$ to the twins. Since $t_e$ has only one neighbor on
$H$ and $t_i$ has only one neighbor other than $b_3$, $H$ begins with
\(
t_eua_1b_1a_2vt_ib_3,
\)
up to interchanging the roles of $a_1$ and $a_2$; see
Figure~\ref{fig:k7-forced}.

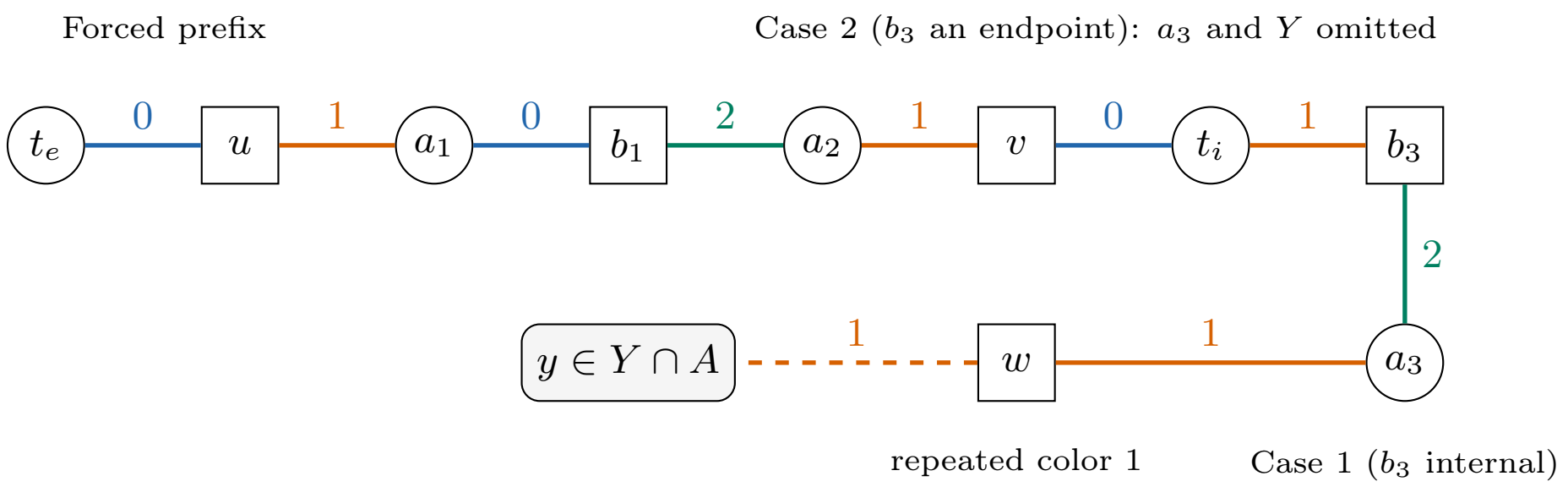
\begin{figure}[htbp]
\centering
\begin{tikzpicture}[
    x=0.95cm,y=1cm,
    every node/.style={font=\small},
    av/.style={circle,draw,fill=white,minimum size=6mm,inner sep=1pt},
    bv/.style={rectangle,draw,fill=white,minimum size=6mm,inner sep=1pt},
    pe/.style={line width=1pt}
]
\node[anchor=west,font=\scriptsize]
    at (0,2.0) {Forced prefix};
\node[anchor=east,font=\scriptsize]
    at (11.6,2.0) {Case~2 ($b_3$ an endpoint): $a_3$ and $Y$ omitted};

\node[av] (te) at (0,1.1) {$t_e$};
\node[bv] (u) at (1.6,1.1) {$u$};
\node[av] (a1) at (3.2,1.1) {$a_1$};
\node[bv] (b1) at (4.8,1.1) {$b_1$};
\node[av] (a2) at (6.4,1.1) {$a_2$};
\node[bv] (v) at (8,1.1) {$v$};
\node[av] (ti) at (9.6,1.1) {$t_i$};
\node[bv] (b3) at (11.2,1.1) {$b_3$};

\draw[pe,edgezero] (te)--node[above]{$0$}(u);
\draw[pe,edgeone] (u)--node[above]{$1$}(a1);
\draw[pe,edgezero] (a1)--node[above]{$0$}(b1);
\draw[pe,edgetwo] (b1)--node[above]{$2$}(a2);
\draw[pe,edgeone] (a2)--node[above]{$1$}(v);
\draw[pe,edgezero] (v)--node[above]{$0$}(ti);
\draw[pe,edgeone] (ti)--node[above]{$1$}(b3);

\node[av] (a3) at (11.2,-.6) {$a_3$};
\node[bv] (w) at (8,-.6) {$w$};
\node[draw,rounded corners,fill=gray!8,minimum height=6mm]
    (y) at (4.8,-.6) {$y\in Y\cap A$};

\draw[pe,edgetwo] (b3)--node[right]{$2$}(a3);
\draw[pe,edgeone] (a3)--node[above]{$1$}(w);
\draw[pe,edgeone,dashed] (w)--node[above]{$1$}(y);

\node[anchor=north,font=\scriptsize]
    at (8,-1.15) {repeated color $1$};
\node[anchor=north,font=\scriptsize]
    at (11.2,-1.15) {Case~1 ($b_3$ internal)};

\end{tikzpicture}
\caption{The forced prefix and its only possible continuation.}
\label{fig:k7-forced}
\end{figure}

As a result,
in Case~2, this prefix already joins the two endpoints of $H$ while
omitting $a_3$ and all vertices of $Y$, a contradiction. In Case~1,
the path continues from $b_3$ to $a_3$. Since the terminal $Y$-block
begins at a vertex of $Y\cap A$, the next vertex is not in $Y$, so it
is the remaining vertex $w\in W\setminus\{u,v\}$. If $w=b_0$, the
color-$2$ edge $a_3b_0$ repeats the color of $b_3a_3$. Hence
$w\in\{b_2,b_4\}$ and $c(a_3w)=1$. All vertices outside $Y$ have now
been visited, so the next edge joins $w$ to a vertex of $Y\cap A$,
again with color $1$. This contradicts proper coloring at $w$.

We now consider the general case.
For $n\geq7$, put $r=n-5\geq2$, and replace $Y$ by
$Y=\{s_1,\ldots,s_r,t_1,\ldots,t_r\}$, where $s_i\in A$ and
$t_j\in B$. Keep the displayed coloring on the core
$\{a_0,\ldots,a_4,b_0,\ldots,b_4\}$. Set $c(s_it_j)=0$ if $i=j$ and
$c(s_it_j)=1$ otherwise. Every edge from $s_i$ to a core vertex has
color $1$, and every edge from a core vertex to $t_j$ has color $0$,
except $c(a_3t_2)=2$. For $n=7$, this is the coloring above, with
$s_1,s_2,t_1,t_2$ corresponding to $a_5,a_6,b_5,b_6$. Since $r\geq2$,
the cycle $s_1t_1s_2t_2\cdots s_rt_rs_1$ alternates between colors $0$
and $1$, and together with $P$ and $C_1$ it forms a PC $1$-path-cycle
factor. Every row and every column contains at least two colors, and
each $s_i$ has exactly two; hence $\delta^c(G)=2$. The preceding
nonexistence argument applies unchanged: $Y$ is balanced, the labels
$\ell(s_i)=1$ and $\ell(t_j)=0$ satisfy the same conditions, and
$a_3t_2$ is the unique edge between $Y$ and $V(G)\setminus Y$ whose
color differs from the label of its endpoint in $Y$. The
remaining facts used in the argument involve only entries of the core
and the new entries $c(s_ib)=1$ and $c(at_j)\in\{0,2\}$, where $a$ and
$b$ are core vertices; hence they persist: the rows of the twins $a_0$
and $a_4$ coincide; the only entry of the column of $b_3$ different from
$1$ is $c(a_3b_3)=2$, and in the columns of $W$ the entries different
from $1$ are $0$ in the rows of the twins and $c(a_3b_0)=2$. Both cases
therefore force the same prefix and lead to the same contradiction.
\end{proof}

Together with the two-color result~\cite{BYTwo},
Theorem~\ref{thm:counterexample} shows that three is the minimum
number of colors for which Conjecture~\ref{conj:bipartite-BJG} fails,
even under the additional assumption $\delta^c(G)\geq2$; for
one color, a PC subgraph is a matching, so a PC $1$-path-cycle factor
exists only for $n=1$, where it is itself a PCHP. The
construction contains the monochromatic cycle $a_0b_0a_4b_2a_0$, and
hence lies outside the forbidden-subgraph hypotheses of~\cite{BYC4}.
It also lies outside the class $\mathcal F$: the monochromatic color
component of color $0$ containing $a_0$ contains $a_1$, $b_0$, and
$b_2$, but both $a_1b_0$ and $a_1b_2$ have color $1$, so this
component is neither a path, nor an even cycle, nor a complete bipartite
graph with a matching removed. These observations motivate the search
for structural restrictions on the coloring under which a PC
$1$-path-cycle factor guarantees a PCHP.

\section{A Characterization via Monochromatic Color Components}
\label{sec:monochromatic-components}

The counterexample in Section~\ref{sec:counterexamples} shows that a
PC $1$-path-cycle factor need not yield a PCHP, even with minimum
color-degree two. In this section, we identify a structural condition
on the monochromatic color components under which the
characterization in Conjecture \ref{conj:bipartite-BJG} holds.

An edge-coloring of a bipartite graph is called a
\emph{bi-equivalence coloring} if every monochromatic color component
is complete bipartite. Such colorings have been studied in connection
with monochromatic coverings and
bicliques~\cite{ChenEtAl2012,DeBiasioEtAl2020}. We consider the
following broader class.

\begin{definition}\label{def:class-F}
\emph{Let $\mathcal F$ denote the class of edge-colored bipartite graphs in
which every monochromatic color component is a path, an even cycle,
or a complete bipartite graph with a possibly empty matching removed.}
\end{definition}
\noindent
Since a complete bipartite graph is obtained by removing the empty
matching, $\mathcal F$ contains all bi-equivalence colorings. In
addition, it allows monochromatic color components that are paths,
cycles, or complete bipartite graphs with a nonempty matching removed, such as
crown graphs. (The same three types of connected bipartite graphs appear
in the characterization of fork-free bipartite
graphs~\cite{Alekseev2004}.)

Throughout this section, $G=K_{m,n}$ is an edge-colored complete
bipartite graph in $\mathcal F$. Since $G$ is complete bipartite, any
two vertices in opposite partite sets are adjacent; in particular, a
cut path of a cycle can be joined to any vertex lying in the partite
set opposite to its corresponding endpoint. The first lemma shows
that, when a vertex of a path sends only edges of one color to a cycle
and this color also occurs on the cycle, the structure of
$\mathcal F$ allows the cycle to be inserted into the path.

\begin{lemma}\label{lem:uniform-completion}
Let $G=K_{m,n}\in\mathcal F$, and let $P$ be a PC path and $C$ a PC
cycle in $G$ that are vertex-disjoint. Suppose that $xy\in E(P)$ has
color $\alpha$ and that every edge between $x$ and $V(C)$ has color
$\alpha$. If $\alpha\in c(C)$, then $C$ can be inserted into $P$ in
place of $xy$, using two joining edges of color $\alpha$, to obtain a
PC path on $V(P)\cup V(C)$.
\end{lemma}

\begin{proof}
Since $\alpha\in c(C)$, there is an edge $uv\in E(C)$ of color
$\alpha$; we may assume that $u$ lies in the same partite set as $x$.
Let $w$ be the other neighbor of $u$ on $C$. Since $C$ is PC,
$c(uw)\ne\alpha$. The vertices $v$, $w$, and $y$ lie in the partite
set opposite to $x$, and the edges $xv$, $xw$, $xy$, and $uv$ all have
color $\alpha$. Hence $x,y,u,v,w$ belong to the same monochromatic
color component $H$ of color $\alpha$. Since $x$ has at least three
neighbors in $H$, the component $H$ is neither a path nor a cycle.
Therefore $H$ is a complete bipartite graph with a matching removed.
The vertices $u$ and $w$ lie in opposite partite sets of $H$, and
$uw\notin E(H)$. Thus $uw$ is the removed matching edge at $u$, and
every other vertex of $H$ in the partite set opposite to $u$ is
adjacent to $u$ in $H$. In particular, since $y\ne w$, we have
$c(uy)=\alpha$. 

Let $Q$ be the cut path of $C$ at $uv$, which runs from $v$ to $u$.
Replace the edge $xy$ of $P$ by the segment $x\,Q\,y$, that is, by
the joining edge $xv$, the path $Q$, and the joining edge $uy$. The
result is a path on $V(P)\cup V(C)$, and it is PC. Both joining edges
have color $\alpha$. The edges of $Q$ at $v$ and at $u$ are cycle
edges adjacent to $uv$, so their colors differ from
$c(uv)=\alpha$ because $C$ is PC. Likewise, the edges of $P$ at $x$
and at $y$ other than $xy$, when present, have colors different from
$c(xy)=\alpha$ because $P$ is PC. Hence the two edges at each of
$x$, $v$, $u$, and $y$ have distinct colors.
\end{proof}

The next lemma is the key merging step. Its proof propagates the
failure of all attachments and insertions along the path, and
Lemma~\ref{lem:uniform-completion} prevents the propagated colors
from occurring on the cycle. This propagation argument
extends the one used for two colors
in~\cite{BYTwo}; the new ingredient is
that the structure of $\mathcal F$ forces the propagated colors to be
absent from the cycle.

\begin{lemma}\label{lem:component-merge}
Let $G=K_{m,n}\in\mathcal F$, and let $P$ be a PC path and $C$ a PC
cycle in $G$ that are vertex-disjoint. Then $C$ can be attached to or
inserted into $P$ so that the result is a PC path on
$V(P)\cup V(C)$.
\end{lemma}

\begin{proof}
Write $P=x_0x_1\cdots x_{q-1}$.
For the trivial case $q=1$,
choose a vertex $v\in V(C)$ in the partite set opposite to $x_0$.
Since $C$ is PC, at least one of the two cycle edges at $v$ has color
different from $c(x_0v)$. Let $e$ be the other cycle edge at $v$. The
cut path of $C$ at $e$ ends at $v$, and attaching it to $x_0$ by the
joining edge $x_0v$ gives a PC path.
For nontrivial cases $q\ge2$,
write $c_i=c(x_ix_{i+1})$ for $0\le i\le q-2$. Suppose, for a
contradiction, that neither an attachment nor an insertion yields a
PC path. We first observe the following.

\begin{equation}\label{eq:cut-choice}\tag{$\ast$}
\parbox{0.82\textwidth}{\emph{For every $v\in V(C)$ and every vertex
$x$ adjacent to $v$, there is a cut path of $C$ ending at $v$ whose
edge at $v$ has a color different from $c(xv)$.}}
\end{equation}
\vskip 0.3cm
\noindent
Indeed, the two cycle edges at $v$ have distinct colors, so at least
one of them differs from $c(xv)$; deleting the other one gives the
required cut path.

We now prove by induction on $i$ that, for every $0\le i\le q-2$,
\begin{equation}\label{eq:propagation}
\text{every edge between $x_i$ and $V(C)$ has color $c_i$, and }
c_i\notin c(C).
\end{equation}
\noindent
For $i=0$, suppose that some edge $x_0v$ has color different from
$c_0$. By $(\ast)$, attaching a suitable cut path of $C$ to $x_0$ by
the joining edge $x_0v$ gives a PC path, a contradiction. Hence every
edge between $x_0$ and $V(C)$ has color $c_0$, and
Lemma~\ref{lem:uniform-completion}, applied to $x_0x_1$, shows that
$c_0\notin c(C)$, since otherwise an insertion would be possible.

Now let $1\le i\le q-2$, and assume that~\eqref{eq:propagation} holds
for $i-1$. Suppose that some edge $x_iv$ has color different from
$c_i$, and let $Q$ be the cut path given by $(\ast)$, running from
$v'$ to $v$. The vertex $v$ lies in the partite set of $x_{i-1}$, and
$v'$, being the other endpoint of a path of even order, lies in the
partite set of $x_i$. Insert $C$ in place of $x_{i-1}x_i$ by the
segment $x_{i-1}\,Q\,x_i$. We check proper coloring at the four
affected vertices:
\begin{itemize}
\item at $x_{i-1}$: the joining edge $x_{i-1}v'$ has color $c_{i-1}$
by the induction hypothesis, while the preceding path edge, if
present, has color $c_{i-2}\ne c_{i-1}$;
\item at $v'$: the joining edge has color $c_{i-1}\notin c(C)$, so it
differs from the edge of $Q$ at $v'$;
\item at $v$: the edge of $Q$ at $v$ differs from $c(x_iv)$ by the
choice of $Q$;
\item at $x_i$: the joining edge has color $c(x_iv)\ne c_i$, the color
of the following path edge $x_ix_{i+1}$.
\end{itemize}
This gives a PC insertion, a contradiction. Hence every edge between
$x_i$ and $V(C)$ has color $c_i$, and
Lemma~\ref{lem:uniform-completion}, applied to $x_ix_{i+1}$, yields
$c_i\notin c(C)$. This completes the induction.

Finally, consider the endpoint $x_{q-1}$. If some edge $x_{q-1}v$ had
color different from $c_{q-2}$, then attaching a cut path of $C$ ending
at $v$, chosen as in $(\ast)$, would give a PC path. Hence every edge
between $x_{q-1}$ and $V(C)$ has color $c_{q-2}$. Let $Q$ be any cut
path of $C$, oriented so that it runs from a vertex $v'$ in the
partite set of $x_{q-1}$ to a vertex $v$ in the partite set of
$x_{q-2}$. Inserting $C$ in place of $x_{q-2}x_{q-1}$ by the segment
$x_{q-2}\,Q\,x_{q-1}$ gives a PC path: both joining edges have
color $c_{q-2}\notin c(C)$ by~\eqref{eq:propagation}, the preceding
path edge at $x_{q-2}$, if present, has color $c_{q-3}\ne c_{q-2}$,
and $x_{q-1}$ becomes an endpoint. This contradiction completes the
proof.
\end{proof}

We record how the endpoints change in Lemma~\ref{lem:component-merge}.
An insertion preserves both endpoints of $P$. An attachment replaces
an endpoint of $P$ by an endpoint of the cut path; since the cut path
has even order, the new endpoint lies in the same partite set as the
old one. In particular, when $q=1$, both endpoints of the resulting
path lie in the partite set of $x_0$. Thus the merging step applies
to paths of either parity and preserves the difference
$|V(P)\cap A|-|V(P)\cap B|$.

\begin{theorem}\label{thm:component-characterization}
For all $m,n\geq1$, an edge-colored complete bipartite graph
$K_{m,n}\in\mathcal F$ contains a PCHP if and only if it admits a PC
$1$-path-cycle factor. Moreover, any such factor can be transformed
into a PCHP by absorbing its cycles one at a time, in an arbitrary
order.
\end{theorem}

\begin{proof}
Necessity is immediate, since a PCHP is itself a PC $1$-path-cycle
factor. Conversely, let $F=\{P,C_1,\ldots,C_k\}$ be a PC
$1$-path-cycle factor of $G=K_{m,n}\in\mathcal F$, and set $Q_0=P$.
For $i=1,\ldots,k$, the path $Q_{i-1}$ and the cycle $C_i$ are
vertex-disjoint, so Lemma~\ref{lem:component-merge} yields a PC path
$Q_i$ on $V(Q_{i-1})\cup V(C_i)$. After $k$ steps, $Q_k$ spans $V(G)$
and is therefore a PCHP. The cycles may be processed in any order.
Since every cycle is balanced, each $Q_i$ satisfies
$|V(Q_i)\cap A|-|V(Q_i)\cap B|=m-n$. In particular, if $|m-n|=1$,
both endpoints of the final path lie in the larger partite set.
\end{proof}

Since every bi-equivalence coloring belongs to $\mathcal F$, we
obtain the following consequence.

\begin{corollary}\label{cor:biequivalence}
For all $m,n\geq1$, a complete bipartite graph $K_{m,n}$ with a
bi-equivalence coloring contains a PCHP if and only if it admits a PC
$1$-path-cycle factor.
\end{corollary}

We conclude this section with two examples illustrating the scope of
Theorem~\ref{thm:component-characterization}. The first shows that
the characterization applies beyond the earlier restrictions on the
number of colors and on monochromatic $4$-cycles, and beyond
bi-equivalence colorings. The second illustrates the nearly balanced
case $|m-n|=1$.

\begin{example}\label{ex:scope-characterization}
Consider the three-edge-coloring of $K_{3,3}$ with partite sets
$A=\{a_0,a_1,a_2\}$ and $B=\{b_0,b_1,b_2\}$ given by the color matrix
\[
\begin{pmatrix}
0&0&0\\
0&0&1\\
2&1&2
\end{pmatrix}.
\]
The monochromatic color component of color $0$ is a $K_{2,3}$ with
the edge $a_1b_2$ removed, the two monochromatic color components of
color $1$ are single edges, and the monochromatic color component of
color $2$ is a $P_3$. Thus the coloring belongs to $\mathcal F$ but is
not a bi-equivalence coloring. Moreover, the graph contains the
monochromatic cycle $a_0b_0a_1b_1a_0$ and the semi-monochromatic
cycle $a_0b_0a_2b_2a_0$, with color sequence $(0,2,2,0)$. Since the
coloring uses three colors and contains a monochromatic $C_4$,
neither the two-color result~\cite{BYTwo} nor the forbidden-$C_4$
result~\cite{BYC4} applies.

The path $P=a_0b_0$ and the cycle $C=a_1b_1a_2b_2a_1$ form a PC
$1$-path-cycle factor. Every edge between $a_0$ and $V(C)$ has
color~$0$, which also appears on $C$; this is precisely the situation
of Lemma~\ref{lem:uniform-completion}. Inserting the cut path of $C$
at $a_1b_1$ in place of $a_0b_0$ yields the PCHP
\(
a_0b_1a_2b_2a_1b_0,
\)
with color sequence $(0,1,2,1,0)$.
\end{example}

\begin{example}\label{ex:unbalanced}
Consider $K_{4,3}$ with partite sets $A=\{a_0,a_1,a_2,a_3\}$ and
$B=\{b_0,b_1,b_2\}$, and the color matrix
\[
\begin{pmatrix}
0&0&3\\
0&1&3\\
2&0&1\\
3&2&0
\end{pmatrix}.
\]
Every monochromatic color component is a path, so the coloring
belongs to $\mathcal F$. The path $P=a_0b_2a_2$ and the cycle
$C=a_1b_0a_3b_1a_1$ form a PC $1$-path-cycle factor, with color
sequences $(3,1)$ and $(0,3,2,1)$, respectively. Inserting the cut
path of $C$ at $a_3b_1$ in place of $a_0b_2$ yields the PCHP
\(
a_0b_1a_1b_0a_3b_2a_2,
\)
with color sequence $(0,1,0,3,0,1)$. Both endpoints lie in the larger
partite set $A$. Transposing the matrix gives the corresponding
example for $K_{3,4}$.
\end{example}

\section{Algorithmic Consequences}\label{sec:algorithms}

We now turn the constructive proof of
Theorem~\ref{thm:component-characterization} into efficient
algorithms. We show that membership in $\mathcal F$ can be recognized
in $O(mn)$ time and that, given a PC $1$-path-cycle factor of
$K_{m,n}\in\mathcal F$, a PCHP can be constructed within the same
time bound, which is linear in the size of the color
matrix. Combined with a polynomial-time algorithm for finding a
maximum PC $1$-path-cycle subgraph \cite{BJG}, this yields a
polynomial-time algorithm for deciding the existence of a PCHP in
$\mathcal F$.

Throughout this section, the edge-coloring of $K_{m,n}$ is given by
its $m\times n$ color matrix, with constant-time access to its
entries. We assume that the colors are given as integers in $\{1,\ldots,s\}$, where
$s=|c(E(G))|\le mn$, and each color occupies one word, so color
comparisons take constant time. Paths and cycles are represented by
ordered vertex lists, and running times are measured in word
operations.

\subsection{Recognizing the class \texorpdfstring{$\mathcal F$}{F}}

Let $H$ be a monochromatic color component with partite sets $A_H$
and $B_H$. If $\Delta(H)\le2$, then $H$ is a path or a cycle, since
it is connected, and every cycle in a bipartite graph is even.
Otherwise, $H$ is a complete bipartite graph with a matching removed
if and only if every vertex of $H$ misses at most one vertex of the
opposite partite set of $H$, that is,
\begin{equation}\label{eq:degree-condition}
d_H(v)\ge|B_H|-1 \quad (v\in A_H)
\qquad\text{and}\qquad
d_H(v)\ge|A_H|-1 \quad (v\in B_H).
\end{equation}
Indeed, \eqref{eq:degree-condition} states exactly that the non-edges
of $H$ between $A_H$ and $B_H$ form a matching.

\begin{algorithm}
\caption{\textsc{RecognizeClass}$(G)$}
\label{alg:recognize}
\begin{algorithmic}[1]
\Require The color matrix of $G=K_{m,n}$
\Ensure \textbf{true} if $G\in\mathcal F$, and \textbf{false} otherwise
\State Group the edges by color
\For{each color $i$}
    \State Build adjacency lists of $G_i$ on the vertices incident with color $i$
    \For{each monochromatic color component $H$ of $G_i$}
        \State Compute $A_H$, $B_H$, and $d_H(v)$ for all $v\in V(H)$
        \If{$\Delta(H)>2$ and~\eqref{eq:degree-condition} fails}
            \State \Return \textbf{false}
        \EndIf
    \EndFor
\EndFor
\State \Return \textbf{true}
\end{algorithmic}
\end{algorithm}

\begin{proposition}\label{prop:recognition}
Algorithm~\ref{alg:recognize} recognizes membership in $\mathcal F$
in $O(mn)$ time.
\end{proposition}

\begin{proof}
Correctness follows from the characterization above. Grouping the
edges by color, for instance by counting sort, takes $O(mn+s)$ time.
For each color $i$, the adjacency lists of $G_i$ are built only on the
vertices incident with color $i$, using timestamped arrays, so that
no array of size $m+n$ is reinitialized for each color. Since each
edge contributes two incidences, building all lists, identifying all
monochromatic color components by graph search, and checking their
degrees take $O(mn+m+n+s)$ time in total. As $s\le mn$ and
$m+n\le mn+1$ for $m,n\ge1$, the running time is $O(mn)$.
\end{proof}

\subsection{Constructing a PCHP by cycle absorption}

We next implement the merging step of
Lemma~\ref{lem:component-merge}. Let $P=x_0\cdots x_{q-1}$ be a PC
path and $C=v_0\cdots v_{\ell-1}v_0$ a vertex-disjoint PC cycle. For
$0\le j<\ell$ and $\varepsilon\in\{-1,1\}$, let
\(
T(j,\varepsilon)=v_jv_{j+\varepsilon}\cdots v_{j+(\ell-1)\varepsilon},
\)
with indices taken modulo $\ell$. Thus $T(j,\varepsilon)$ is the cut
path of $C$ at $v_jv_{j-\varepsilon}$, traversed from $v_j$, and these
$2\ell$ sequences list every cut path of $C$ in both orientations. For
a position $b\in\{0,\ldots,q\}$, the \emph{candidate}
\(
R(j,\varepsilon,b)=x_0\cdots x_{b-1}\,T(j,\varepsilon)\,x_b\cdots x_{q-1}
\)
is obtained by attaching $C$ to $P$ when $b\in\{0,q\}$, and by
inserting $C$ in place of $x_{b-1}x_b$ when $1\le b\le q-1$; empty
portions are omitted. By Lemma~\ref{lem:component-merge}, at least one
candidate is a PC path.

Since $P$ and $T(j,\varepsilon)$ are PC, a candidate is a PC path if
and only if each joining edge exists, that is, its endpoints lie in
opposite partite sets, and at each endpoint of a joining edge the two
incident edges of the candidate have distinct colors. These conditions
involve at most two joining edges and the four vertices $x_{b-1}$,
$v_j$, $v_{j-\varepsilon}$, and $x_b$, whenever present. Hence they
can be tested in constant time from the color matrix and the indexed
vertex lists, without constructing the candidate. We denote this test
by $\operatorname{Compatible}(P,C,j,\varepsilon,b)$. The procedure
\textsc{Merge} (Algorithm~\ref{alg:merge}) scans the candidates and
returns the first one that passes the test.

\begin{algorithm}[htbp]
\caption{Merging a PC cycle into a PC path}
\label{alg:merge}
\begin{algorithmic}[1]
\Require Vertex-disjoint PC path $P=x_0\cdots x_{q-1}$ and PC cycle
$C=v_0\cdots v_{\ell-1}v_0$ in $G=K_{m,n}\in\mathcal F$
\Ensure A PC path on $V(P)\cup V(C)$
\Function{Merge}{$P,C$}
    \For{$j=0,\ldots,\ell-1$}
        \For{$\varepsilon\in\{-1,1\}$}
            \For{$b=0,\ldots,q$}
                \If{$\operatorname{Compatible}(P,C,j,\varepsilon,b)$}
                    \State \Return $R(j,\varepsilon,b)$
                \EndIf
            \EndFor
        \EndFor
    \EndFor
\EndFunction
\end{algorithmic}
\end{algorithm}

Given a PC $1$-path-cycle factor $F=\{P,C_1,\ldots,C_k\}$, the
procedure \textsc{ConstructPCHP} (Algorithm~\ref{alg:factor}) starts
with $Q=P$ and absorbs the cycles one at a time by calling
\textsc{Merge}. By Lemma~\ref{lem:component-merge}, each call
succeeds, and the cycles may be processed in any order. After each
call, $Q$ is a PC path on the vertices processed so far; after the
last call, it is a PCHP. If $k=0$, the path $P$ is already
Hamiltonian and is returned unchanged.

\begin{theorem}\label{thm:factor-complexity}
Given a PC $1$-path-cycle factor of an edge-colored complete bipartite
graph $K_{m,n}\in\mathcal F$, Algorithm~\ref{alg:factor} constructs a
PCHP in $O(mn)$ time.
\end{theorem}

\begin{proof}
Correctness follows from Lemma~\ref{lem:component-merge}. Consider a
call of \textsc{Merge} with a path of order $q$ and a cycle of order
$\ell$. It tests at most $2\ell(q+1)$ candidates, each in constant
time, and constructs the returned path once, in $O(q+\ell)$ time.
Hence the call takes $O(\ell(q+1))$ time.

Let $N=m+n$, and let $q_i$ and $\ell_i$ denote the orders of the
current path and of the cycle in the $i$th call. Since $q_i\le N$ and
$\sum_{i=1}^{k}\ell_i\le N$, the total running time is
\[
O\Bigl(N+\sum_{i=1}^{k}\ell_i(q_i+1)\Bigr)=O(N^2),
\]
which also covers the case $k=0$. Finally, the existence of a PC
$1$-path-cycle factor implies $|m-n|\le1$ (see
Section~\ref{sec:preliminaries}), so
\[
N^2=(m+n)^2=4mn+(m-n)^2\le4mn+1,
\]
and the running time is $O(mn)$, which is linear in the
size of the color matrix.
\end{proof}

\begin{algorithm}[htbp]
\caption{Constructing a PCHP from a PC $1$-path-cycle factor}
\label{alg:factor}
\begin{algorithmic}[1]
\Require $G=K_{m,n}\in\mathcal F$ and a PC $1$-path-cycle factor
$F=\{P,C_1,\ldots,C_k\}$ of $G$
\Ensure A PCHP of $G$
\Function{ConstructPCHP}{$G,F$}
    \State $Q\gets P$
    \For{$i=1,\ldots,k$}
        \State $Q\gets\Call{Merge}{Q,C_i}$
    \EndFor
    \State \Return $Q$
\EndFunction
\end{algorithmic}
\end{algorithm}
Bang-Jensen and Gutin~\cite{BJG} showed that a maximum
PC $1$-path-cycle subgraph of an edge-colored multigraph, where maximum refers to the number of covered vertices,
can be found in polynomial time. Combining this with
Theorem~\ref{thm:component-characterization} and
Algorithm~\ref{alg:factor}, we obtain the following.

\begin{corollary}\label{cor:polynomial-decision}
For edge-colored complete bipartite graphs in $\mathcal F$, the
existence of a PCHP can be decided in polynomial time. Moreover, a
PCHP can be constructed in polynomial time whenever one exists.
\end{corollary}

\begin{proof}
If $|m-n|>1$, no Hamiltonian path exists. Otherwise, find a maximum
PC $1$-path-cycle subgraph. If it is not spanning, then $G$ has no PC
$1$-path-cycle factor and hence no PCHP, since a PCHP is itself such
a factor. If it is spanning, it is a PC $1$-path-cycle factor, and
Algorithm~\ref{alg:factor} converts it into a PCHP in $O(mn)$ time.
Both steps take polynomial time. When $|m-n|=1$, both endpoints of this PCHP lie in the larger
partite set. 
\end{proof}

\section{Concluding Remarks}\label{sec:scope}

We have shown that the conjectured extension of the
Bang-Jensen-Gutin Theorem to complete balanced bipartite graphs
fails in general, even for three-edge-colorings with minimum
color-degree two. On the positive side, we established the
characterization for the class $\mathcal F$, defined by restricting
the structure of the monochromatic color components. This class
strictly contains the bi-equivalence colorings, and the
characterization holds not only in the balanced case but also in the
nearly balanced case $|m-n|=1$. The proof is constructive: every PC
$1$-path-cycle factor can be transformed into a PCHP by absorbing its
cycles one at a time, in any order. We also obtained $O(mn)$-time
algorithms for recognizing $\mathcal F$ and for constructing a PCHP
from a given factor, which lead to a polynomial-time decision
algorithm for this class.

A main idea of our approach is that each cycle is absorbed into
the current path independently of the others. The following example
shows that this can fail even when a PCHP exists, and that a PCHP may
then still be obtained by first rearranging several cycles.

\begin{example}\label{ex:cycle-first}
Consider the edge-colored $K_{7,7}$ with bipartition
$A=\{a_0,\ldots,a_6\}$ and $B=\{b_0,\ldots,b_6\}$, whose color matrix
is
\[
\begin{pmatrix}
2&2&2&2&2&2&2\\
1&0&1&0&0&0&0\\
1&2&1&0&0&0&0\\
1&1&1&0&1&2&2\\
1&1&1&1&0&2&2\\
1&1&1&2&2&0&1\\
1&1&1&2&2&1&0
\end{pmatrix}.
\]
The graph admits the PC $1$-path-cycle factor $F=\{P,C_1,C_2\}$, where
\(P=a_0b_0a_1b_1a_2b_2\), \(C_1=a_3b_3a_4b_4a_3\), and
\(C_2=a_5b_5a_6b_6a_5\).

\textbf{First, we show that no single cycle can be absorbed:}
We show that neither $G[V(P)\cup V(C_1)]$ nor $G[V(P)\cup V(C_2)]$
contains a PCHP. Fix $i\in\{1,2\}$, let $Y=V(C_i)$ and
$X'=V(P)\setminus\{a_0\}$, and suppose that $H$ is a PCHP of
$G[V(P)\cup Y]$. Since every edge at $a_0$ has color $2$, the vertex
$a_0$ is an endpoint of $H$. Assign the label $1$ to the vertices of
$Y\cap A$ and the label $0$ to those of $Y\cap B$. Then every edge of
$G[Y]$ has the label of one of its endpoints, and every edge between
$X'$ and $Y$ has the label of its endpoint in $Y$. By
Lemma~\ref{lem:labels}, no $Y$-block of $H$ lies between two vertices
of $X'$. Since the only vertex outside $X'\cup Y$ is the endpoint
$a_0$, the five vertices of $X'$ form a segment $Q$ of $H$, and both
endpoints of $Q$ lie in $B$.

Within $X'$, the only neighbor of $a_1$, and likewise of $a_2$, joined
to it by an edge of color different from $1$ is $b_1$. As $a_1$ and
$a_2$ are internal vertices of $Q$, proper coloring forces $Q$ to
contain the segment $a_1b_1a_2$, up to reversal, and hence both
terminal edges of $Q$ have color $1$. On $H$, the endpoints of $Q$
can be adjacent only to $a_0$ or to vertices of $Y\cap A$, and $a_0$
has only one neighbor on $H$. Since $Y\ne\emptyset$, some endpoint of
$Q$ is adjacent along $H$ to a vertex of $Y\cap A$. This edge has
color $1$, the same as the adjacent terminal edge of $Q$, a
contradiction. Thus neither cycle can be absorbed into $P$, even if
arbitrary edges of the corresponding induced subgraph are allowed.

In particular, for each $i\in\{1,2\}$, the subgraph
$G[V(P)\cup V(C_i)]$ is a three-edge-colored $K_{5,5}$ that admits the
PC $1$-path-cycle factor $\{P,C_i\}$ but has no PCHP. These smaller
examples do not satisfy the minimum color-degree condition of
Theorem~\ref{thm:counterexample}, since all edges at $a_0$ have
color~$2$ and hence $d^c(a_0)=1$.

\textbf{We now rearrange the cycles first and then absorb:}
Deleting $a_3b_3$ and $a_6b_5$ from $C_1\cup C_2$ and adding the
color-$2$ edges $a_3b_5$ and $a_6b_3$ produces the PC cycle
\(
D=b_3a_4b_4a_3b_5a_5b_6a_6b_3,
\)
with color sequence $(1,0,1,2,0,1,0,2)$. Inserting the cut path of
$D$ at $a_6b_3$ into $P$ in place of $a_2b_2$ yields
\(
H=a_0b_0a_1b_1a_2\,b_3a_4b_4a_3b_5a_5b_6a_6\,b_2,
\)
with color sequence $(2,1,0,2,0,1,0,1,2,0,1,0,1)$. Hence $H$ is a
PCHP of $G$.

Note that this coloring lies outside $\mathcal F$. The monochromatic color
component of color $2$ contains $a_0$, which has degree seven in it,
so this component is neither a path nor a cycle. It also contains
$a_2$, $b_0$, and $b_3$, but neither $a_2b_0$ nor $a_2b_3$ has
color~$2$; hence it is not a complete bipartite graph with a matching
removed. Moreover, the graph contains the monochromatic cycle
$a_1b_0a_2b_2a_1$ and the semi-monochromatic cycle
$b_1a_0b_3a_1b_1$, with color sequence $(2,2,0,0)$. Hence the example
lies outside the hypotheses of~\cite{BYTwo,BYC4} and of
Theorem~\ref{thm:component-characterization}.
\end{example}

Therefore, Example~\ref{ex:cycle-first} shows that the failure of individual
cycle absorption does not prevent the existence of a PCHP:
rearranging several components simultaneously may succeed where
successive path-cycle mergers fail. This suggests that the structural
conditions ensuring a PCHP extend beyond those that support
independent cycle absorption, and leads to the following question.

\begin{question}\label{ques:general-construction}
\emph{Under what structural conditions on an edge-colored complete balanced
bipartite graph does the existence of a PC $1$-path-cycle factor
guarantee the existence of a PCHP?}
\end{question}

For complete graphs, Bang-Jensen, Gutin and
Yeo~\cite{BangJensenGutinYeo} showed,
before the BJG Theorem was established, that a PC cycle factor, that
is, a
spanning collection of vertex-disjoint PC cycles, already
guarantees a PCHP. Deleting an edge of one cycle of a PC cycle factor
yields a PC $1$-path-cycle factor, so this statement is weaker than the
BJG Theorem. Our counterexamples do not settle its bipartite analogue,
since they admit no PC cycle factor: in such a factor, both twins $a_0$
and $a_4$ of Theorem~\ref{thm:counterexample} would use their unique
color-$1$ edges to $b_3$, and in the $K_{5,5}$ examples of
Example~\ref{ex:cycle-first} the vertex $a_0$ is incident only with
edges of color $2$. For two colors, the answer is positive
by~\cite{BYTwo}. PC cycle factors in edge-colored balanced bipartite
graphs were recently studied in~\cite{HanEtAl}.

\begin{question}\label{ques:cycle-factor}
\emph{Does every edge-colored complete balanced bipartite graph with a PC
cycle factor contain a PCHP?}
\end{question}

Since Conjecture~\ref{conj:bipartite-BJG} fails, the existence of a
PC $1$-path-cycle factor, which can be tested in polynomial
time~\cite{BJG}, no longer decides the existence of a PCHP in
edge-colored complete balanced bipartite graphs. To the
best of our knowledge, the computational complexity of deciding whether
an edge-colored $K_{n,n}$ has a PCHP is open; it is polynomial for two
colors~\cite{BYTwo} and for colorings in $\mathcal F$
(Corollary~\ref{cor:polynomial-decision}).

A related direction is to develop constructive methods that rearrange
several components of a PC $1$-path-cycle factor simultaneously.
Since a graph may admit many such factors, it would also be
interesting to study transformations between factors and criteria for
selecting a factor from which a PCHP can be constructed.

\bibliographystyle{elsarticle-num}
\bibliography{reference}
\end{document}